\documentclass[12pt, final]{l4dc2022}
\renewcommand{\thefootnote}{*}
\newtheorem{assum}{Assumption}
\usepackage{times}

\title[Data-driven Control with Quantized Feedback]{Data-driven Control of Unknown Linear Systems via Quantized Feedback}

\author{\Name{Feiran Zhao}\footnotemark[1] \Email{zhaofr18@mails.tsinghua.edu.cn}\\
	\Name{Xingchen Li}\footnotemark[1] \Email{lixc21@mails.tsinghua.edu.cn}\\
	\Name{Keyou You} \Email{youky@tsinghua.edu.cn}\\
	\addr Department of Automation and BNRist, Tsinghua University}

\begin{document}
\footnotetext[1]{Equal contribution.}
\maketitle

\begin{abstract}
Control using quantized feedback is a fundamental approach to system synthesis with limited communication capacity. In this paper, we address the stabilization problem for unknown linear systems with logarithmically quantized feedback, via a direct data-driven control method. By leveraging a recently developed matrix S-lemma, we prove a sufficient and necessary condition for the existence of a common stabilizing controller for all possible dynamics consistent with data, in the form of a linear matrix inequality. Moreover, we formulate semi-definite programming to solve the coarsest quantization density. By establishing its connections to unstable eigenvalues of the state matrix, we further prove a necessary rank condition on the data for quantized feedback stabilization. Finally, we validate our theoretical results by numerical examples.
\end{abstract}


\section{Introduction}

Motivated by limited communication capacity in physical systems, control with quantized feedback \cite{elia2001stabilization,fu2005sector,you2011attainability,kang2015coarsest,zhou2019adaptive} has been an active research area for decades, dating back to an early work \citep{kalman1956nonlinear} by Kalman. By loosely compressing the system input and output into sectionalized levels through a quantizer, the dynamical system can be stabilized, possibly with $\mathcal{H}_2$ and $\mathcal{H}_{\infty}$ performance guarantees under low communication load. The main line of the existing literature is to understand and mitigate the side effects of quantization. However, an explicit dynamical model is required which might be impossible to obtain in practical scenarios.

Recent years have witnessed a renewed surge in data-driven control for unknown linear systems \cite{de2019formulas,de2021low,van2020data,van2020noisy}. Instead of identifying a descriptive model via system identification, the direct data-driven framework computes a control law from sampled system trajectories. A key feature of this approach is that it does not require the well-known persistent excitation condition, hence shedding light on the case where data is insufficient for system identification. By invoking a data-based uncertainty representation, it has been shown in \cite{van2020noisy} that for informative data, a stabilizing controller can be efficiently found via a Linear Matrix Inequality (LMI) \citep{boyd1994linear}. 

In this paper, we take an initial step towards quantized control in the direct data-driven approach and focus on the quantized stabilization problem for unknown noisy linear systems with a single input. In particular, we limit the input quantizer to be logarithmic (i.e., the quantization levels are linear in logarithmic scale) with a linear state feedback input, which has been proved \citep{elia2001stabilization} to be able to attain the coarsest quantization density for stabilization of deterministic linear systems. By leveraging a recently developed matrix S-lemma \citep{van2020noisy}, we prove a sufficient and necessary condition for the existence of a common stabilizing controller for all possible dynamics that reflect the data, in terms of an LMI. Moreover, we propose a Semi-definite Programming (SDP) to solve a uniform lower bound for the quantization density. By establishing its connections to unstable eigenvalues of the state matrix, we further prove a necessary rank condition on the data for quantized feedback stabilization.

\subsection{Related work}

\textbf{Quantized control.} Control using quantized feedback has drawn increasing attention of the control community from the seminal work \citep{elia2001stabilization}. Its major contribution is to show that for linear systems with a single input, the coarsest quantizer is logarithmic, and the associated quantization density can be computed using the unstable eigenvalues of the state matrix. \cite{fu2005sector} interpret the quantized stabilization problem as an $\mathcal{H}_{\infty}$ control problem based on the sector bound method, and extend it to multiple-input multiple-output (MIMO) systems. Inspired by \cite{fu2005sector}, our work solves a minimax $\mathcal{H}_{\infty}$ control problem to find the coarsest density. More pertinent to our work is the quantized control for linear uncertain systems. \cite{hayakawa2009adaptive} proposes a time-varying adaptive control law for asymptotic stabilization of uncertain noiseless linear systems by solving a series of Riccati equations. However, we aim to stabilize a set of systems consistent with the data using a common quantized state feedback controller. \cite{kang2015coarsest} considers the element uncertainty in the state matrix of a controllable canonical system, and addresses the stabilization problem of systems with two blocks of uncertainties (due to the plant and the quantization) as done in our work. In contrast, we consider a more natural form of uncertainty which is reflected by the data. Other works include \cite{gao2008new,fu2009finite,coutinho2010input,shen2017quantized}, and \cite{corradini2008robust,liu2012sector,yu2016adaptive,zhou2019adaptive} for nonlinear systems.

\textbf{Direct data-driven control.} This line of work originates from the \textit{Willems et al.'s fundamental lemma} proposed by \cite{willems2005note}, which states that the dynamical model can be replaced with data from sufficiently excited systems. Motivated by it, \cite{de2019formulas} represents the dynamics using historical trajectories under the Persistent Excitation (PE) condition and formulate an SDP to solve a stabilizing controller for deterministic linear systems. This data-driven framework is then applied to model predictive control to deal with safety constraints \citep{coulson2019data, berberich2020data}. When the data is corrupted with noises, \cite{de2021low} proposes an SDP with robustness to the noise for stabilization in an ad-hoc way. For the case that data is insufficient to satisfy the PE condition, the seminal work \cite{van2020data} establishes sufficient and necessary conditions on the informativity of the data for several fundamental control problems. \cite{van2020noisy, bisoffi2021trade} solves a stabilizing controller via an LMI from noisy data by leveraging a matrix S-lemma. Other data-driven work includes \cite{guo2021data, xu2021data, rotulo2021online, berberich2020combining}. Our work extends the framework in \cite{van2020noisy} to address quantized stabilization problems and answer fundamental questions regarding the coarsest quantization density and the condition on the noisy data.

\section{Problem formulation}
In this paper, we consider the following discrete-time linear time-invariant system 
\begin{equation}\label{equ:sys}
x(k+1) = Ax(k) + Bu(k)+w(k),
\end{equation}
where $x(k) \in \mathbb{R}^n$ denotes the state, $u(k) \in \mathbb{R}$ is the control input and $\{w(k) \in \mathcal{R}^n\}$ is an uncorrelated noise sequence. As an initial attempt, we focus on the setting that $A \in \mathbb{R}^{n \times n}$ is an unknown state matrix but the input matrix $B \in \mathbb{R}^{n}$ is known \textit{a priori}, and $(A,B)$ are stabilizable.
\begin{figure}[t]
	\centering
	\includegraphics[height=45mm]{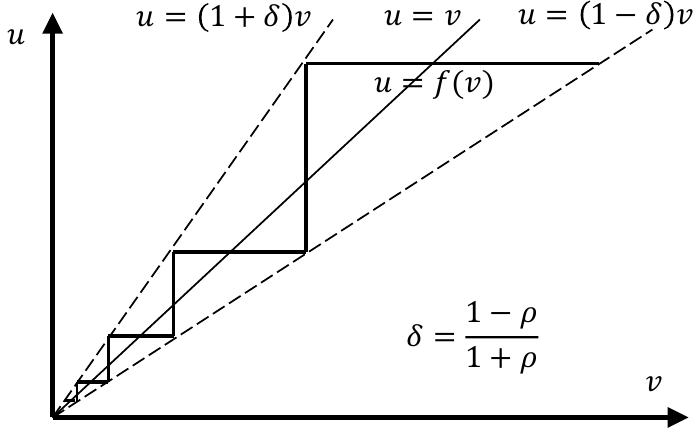}
	\caption{The logarithmic quantizer with density $\rho$.}
	\label{pic:quantizer}
\end{figure}
Since the coarsest quantizer that quadratically stabilizes linear time-invariant systems is logarithmic \citep{elia2001stabilization}, we aim to stabilize (\ref{equ:sys}) via logarithmically quantized state feedback, provided only with a finite length of input and state trajectory $\{x(0), u(0), x(1), u(1), \dots, x(T)\}$. Define the quantization level $u_i = \rho^{-i} u_0, ~i \in \{1,2,\dots\}$ with a density $0<\rho<1$. Then, the quantized controller has the form of 
\begin{equation}\label{def:quantizer}
u(k) = f(v(k)), ~~ v(k) = Kx(k),
\end{equation}
where a logarithmic quantizer $f(\cdot): \mathbb{R} \rightarrow \mathbb{R}$ (see Fig. \ref{pic:quantizer}) is given by
\begin{equation}\label{def:delta}
f(v)= \begin{cases}u_{i}, & \text { if } \frac{1}{1+\delta} u_{i}<v \leq \frac{1}{1-\delta} u_{i}, v>0 \\ 0, & \text { if } v=0 \\ -f(-v), & \text { if } v<0\end{cases} ~~\text{with}~~ \delta = \frac{1-\rho}{1+\rho},
\end{equation}
and $K \in \mathbb{R}^{1 \times n}$ is the feedback gain. We make a standard technical assumption \citep{de2019formulas,van2020noisy} on the noise.
\begin{assum}\label{assum}
	The sequence $W = \{w(0), w(1), \dots, w(T-1)\}$ is unknown, yet satisfies a quadratic bound
	\begin{equation}\label{assumption}
	\begin{bmatrix}
	I \\
	W^{\top}
	\end{bmatrix}^{\top}\begin{bmatrix}
	\Phi_{11} & \Phi_{12} \\
	\Phi_{12}^{\top} & \Phi_{22}
	\end{bmatrix}\begin{bmatrix}
	I \\
	W^{\top}
	\end{bmatrix} \geq 0
	\end{equation}
	with a symmetric matrix $\Phi_{11}$ and a negative definite matrix $\Phi_{22}<0$.
\end{assum}
Assumption \ref{assum} is in the form of a quadratic matrix inequality, which is able to capture important prior knowledge on the system such as energy and sample covariance bounds over a finite horizon. For example, let $\Phi_{12} = 0$ and $\Phi_{22} = - I$, it reduces to $WW^{\top} = \sum_{k=0}^{T-1}w(k)w(k)^{\top} \leq \Phi_{11}$.

A main challenge for quantized stabilization is that the unknown matrix $A$ cannot be deduced uniquely from data, which instead constitutes an uncertainty set. In fact, our approach relies on a data-based representation of the uncertainty~\citep{van2020noisy}. Define the data matrices $X_{-} = \{x(0), x(1) , \dots, x(T-1)\}$, $U = \{u(0), u(1), \dots, u(T-1)\}$, $X_{+} = \{x(1), x(2) , \dots, x(T)\}$, which are constrained by system dynamics
$
X_{+} = AX_{-}+ BU + W.
$
Substituting $W$ into (\ref{assumption}), it follows that all possible $A$ consistent with the data must satisfy a quadratic matrix inequality

\begin{equation}\label{equ:newset}
\begin{bmatrix}
I \\
A^{\top} 
\end{bmatrix}^{\top}
\begin{bmatrix}
I & X_{U} \\
0 & -X_{-} \\
\end{bmatrix}\begin{bmatrix}
\Phi_{11} & \Phi_{12} \\
\Phi_{12}^{\top} & \Phi_{22}
\end{bmatrix}\begin{bmatrix}
I & X_{U} \\
0 & -X_{-} \\
\end{bmatrix}^{\top}
\begin{bmatrix}
I \\
A^{\top} 
\end{bmatrix} \geq 0,
\end{equation}
with $X_{U} = X_{+}-BU$, which defines the uncertainty set 
$\Sigma = \{ A |  (\ref{equ:newset})~\text{holds} \}.$

Our goal is now to design a controller in the form of (\ref{def:quantizer}) to stabilize all systems with $A \in \Sigma$. This naturally raises the following fundamental question: can we characterize the conditions for such a controller to existing? Clearly, the data must be sufficiently informative such that the uncertainty set $\Sigma$ is small enough. Moreover, the quantizer $f(\cdot)$ cannot be too coarse, which implies that the quantization density $\rho$ has a uniform lower bound. 

In this paper, we provide an affirmative answer to the question. By converting the stabilization problem to an $\mathcal{H}_{\infty}$ control problem and applying a recently developed matrix S-lemma \citep{van2020noisy}, we derive a sufficient and necessary condition on the data matrices as well as the quantization density in terms of an LMI. Moreover, we show that the coarsest quantizer and the associated controller can be found by solving a minimax $\mathcal{H}_{\infty}$ norm optimization problem, which can be formulated as an efficient SDP. By relating the minimax problem to the unstable eigenvalues of the state matrix, we further prove an explicit rank condition on the data matrices, which is necessary for quantized stabilization.

\section{Sufficient and necessary condition for quantized stabilization}

In this section, we establish a sufficient and necessary condition for (\ref{equ:sys}) to be stabilizable with logarithmically quantized linear feedback. We first convert the stabilization problem to an $\mathcal{H}_{\infty}$ control problem over the uncertainty set $\Sigma$. 

\begin{figure}[t]
	\centering
	\includegraphics[height=30mm]{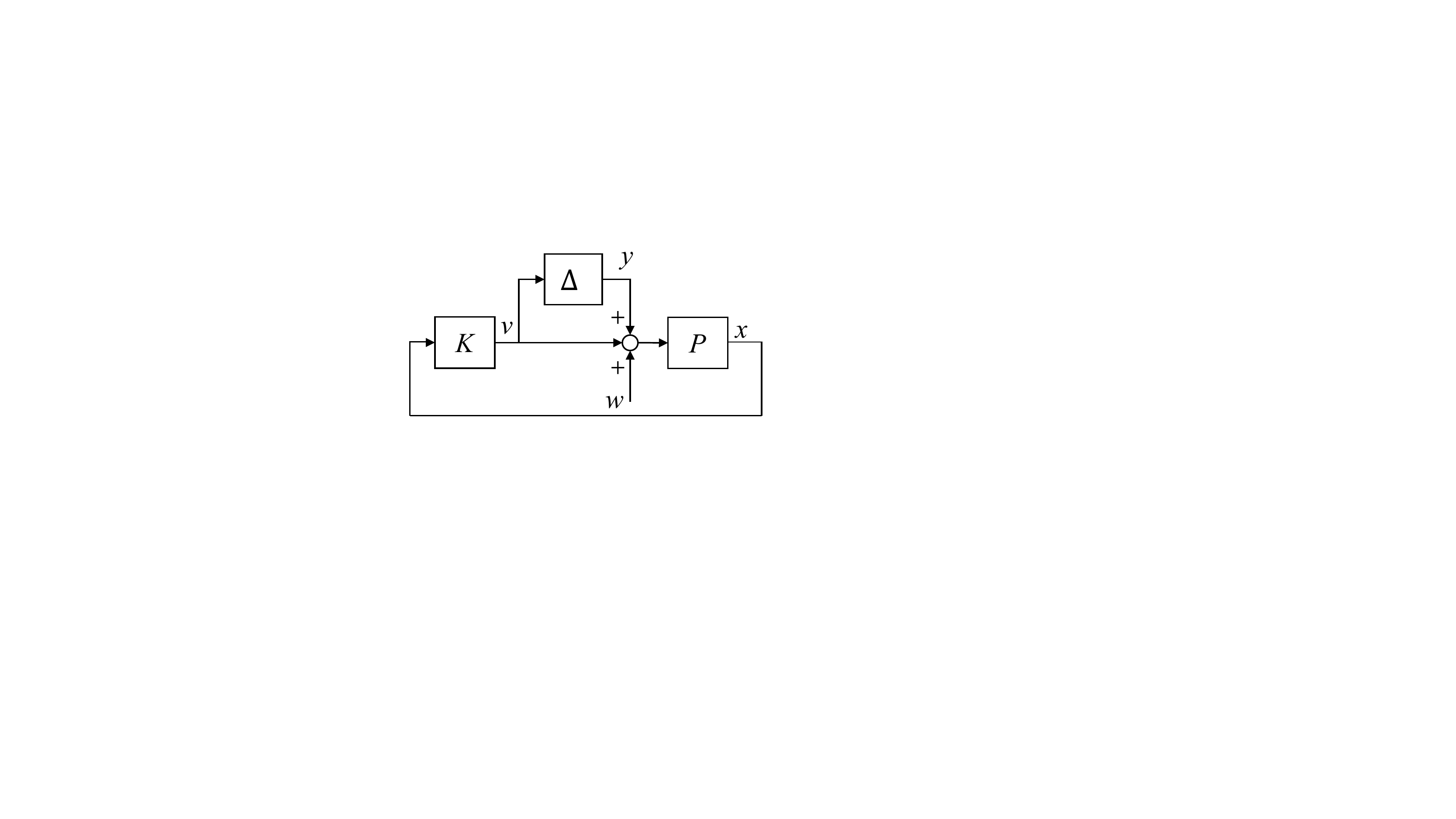}
	\caption{The control diagram of the uncertain system (\ref{equ:uncertain}).}
	\label{pic:diagram}
\end{figure}

It has been shown in \cite{fu2005sector} that (\ref{equ:sys}) is stabilizable via quantized linear state feedback with quantization density $\rho$ if and only if the uncertain system
\begin{equation}\label{equ:uncertain}
x(k+1) = Ax(k) + B(1+\Delta)v(k), ~~\Delta \in [-\delta, \delta]
\end{equation}
is quadratically stabilizable via linear state feedback, where $\delta$ is defined by (\ref{def:delta}). Fig. \ref{pic:diagram} illustrates the control diagram of (\ref{equ:uncertain}) by viewing $\Delta$ as the input uncertainty.
By the well-known small-gain theorem~\citep{zhou1998essentials}, the stabilization of (\ref{equ:uncertain}) is equivalent to an $\mathcal{H}_{\infty}$ control problem. Define the transfer function from signal $y$ to $v$ as
$$G_{A,K}(z) = \left[\begin{array}{c|c}
A +BK & B \\
\hline K & 0
\end{array}\right]=  K(zI-A-BK)^{-1}B.$$
We use $\| G_{A,K}(z)\|_{\infty}$ to denote its $\mathcal{H}_{\infty}$ norm in the Hilbert space. Then, the uncertain system (\ref{equ:uncertain}) is stabilizable via $v = Kx$ if and only if 
\begin{equation}\label{equ:h_inf}
	\delta^{-1} >\| G_{A,K}(z)\|_{\infty}.
\end{equation}

%

Since the condition (\ref{equ:h_inf}) is expressed in the frequency domain, we convert it to an explicit algebraic inequality by applying the bounded real lemma in robust control theory~\citep{zhou1996robust}. 

\begin{lemma}[Bounded real lemma]\label{lem:bounded}
	 $\|G_{A,K}(z)\|_{\infty} < 1/\delta$ if and only if there exists $P>0$ such that 
	\begin{equation}\label{equ:brl}
	I-\delta^2B^{\top}PB >0, ~~\text{and}~~ K^{\top}K + ( A + BK)^{\top}(P^{-1} - \delta^2 BB^{\top})^{-1}( A + BK) < P .
	\end{equation}
\end{lemma}

By Lemma \ref{lem:bounded}, the stabilization problem is now equivalent to solving the $\mathcal{H}_{\infty}$ control problem in (\ref{equ:brl}) subject to $A \in \Sigma$. We show that (\ref{equ:brl}) can be rewritten as a quadratic inequality by a standard change of variables. Let $Y = P^{-1}$, $X = KY$. Pre- and postmultiplying $P^{-1}$ to (\ref{equ:brl}), it holds that
$$
-Y + (AY+BX)^{\top}(Y- \delta^2 BB^{\top})^{-1}(AY+BX) + X^{\top}X < 0.
$$
A Schur complement argument further yields that
$$
\begin{bmatrix}
Y-X^{\top}X & (AY+BX)^{\top} \\
AY+BX & Y-\delta^2 BB^{\top}
\end{bmatrix} > 0,
$$
which again by Schur complement is equivalent to
$$
Y-\delta^2 BB^{\top} - (AY+BX)(Y-X^{\top}X)^{-1} (AY+BX)^{\top} >0, ~~\text{and}~~
Y-X^{\top}X >0.
$$
Reorganizing it in a quadratic form of $A$ with $Z = Y-X^{\top}X>0$, we have that
\begin{equation}\label{equ:performance}
\begin{bmatrix}
I \\
A^{\top}
\end{bmatrix}^{\top}
\begin{bmatrix}
Y- \delta^2 BB^{\top} -BXZ^{-1}X^{\top}B^{\top} & -BXZ^{-1}Y^{\top} \\
-YZ^{-1}X^{\top}B^{\top} & -YZ^{-1}Y^{\top}
\end{bmatrix}
\begin{bmatrix}
I \\
A^{\top}
\end{bmatrix}>0.
\end{equation}   

Then, the quantized stabilization problem leads to the following question: when does (\ref{equ:performance}) hold for all $A$ satisfying (\ref{equ:newset})? Since both (\ref{equ:performance}) and (\ref{equ:newset}) are quadratic inequalities in $[I ~~A^{\top}]$, we apply the matrix-valued S-lemma~\citep{van2020noisy} to derive an LMI condition. 
\begin{lemma}[Matrix S-lemma]\label{lem:matrixS}
	Let $M = \begin{bmatrix}
	M_{11} & M_{12} \\
	M_{12}^{\top} & M_{22}
	\end{bmatrix} \in \mathbb{R}^{(n+n)\times(n+n)}$ and $N = \begin{bmatrix}
	N_{11} & N_{12} \\
	N_{12}^{\top} & N_{22}
	\end{bmatrix} \in \mathbb{R}^{(n+n)\times(n+n)}$ be symmetric matrices with $M_{22}\leq 0$ and $N_{22} \leq 0$. Suppose that $\text{ker}(N_{22}) \subset \text{ker}(N_{12})$ and there exists some matrix $A$ satisfying
	$
		\begin{bmatrix}
		I \\
		A^{\top} 
		\end{bmatrix}^{\top}
		N
		\begin{bmatrix}
		I \\
		A^{\top} 
		\end{bmatrix} > 0
	$ (the so-called generalized Slater condition).
	Then, it follows that $\begin{bmatrix}
	I \\
	A^{\top} 
	\end{bmatrix}^{\top}
	M
	\begin{bmatrix}
	I \\
	A^{\top} 
	\end{bmatrix} > 0$ for all $A$ satisfying $\begin{bmatrix}
	I \\
	A^{\top} 
	\end{bmatrix}^{\top}
	N
	\begin{bmatrix}
	I \\
	A^{\top} 
	\end{bmatrix} \geq 0$ if and only if there exists $\alpha \geq 0$ and $\beta \geq 0$ such that
	$
	M-\alpha N \geq \begin{bmatrix}
	\beta I & 0 \\
	0 & 0
	\end{bmatrix}.
	$
\end{lemma}

In our setting, we define the partitioned matrices
\begin{align}
	M & =\begin{bmatrix}
	M_{11} & M_{12} \\
	M_{12}^{\top} & M_{22}
	\end{bmatrix}= \begin{bmatrix}
	Y- \delta^2 BB^{\top} -BXZ^{-1}X^{\top}B^{\top} & -BXZ^{-1}Y^{\top} \\
	-YZ^{-1}X^{\top}B^{\top} & -YZ^{-1}Y^{\top}
	\end{bmatrix}, \\
	N & = \begin{bmatrix}
	N_{11} & N_{12} \\
	N_{12}^{\top} & N_{22}
	\end{bmatrix} = 
	\begin{bmatrix}
	I & X_{+}-BU \\
	0 & -X_{-} \\
	\end{bmatrix}\begin{bmatrix}
	\Phi_{11} & \Phi_{12} \\
	\Phi_{12}^{\top} & \Phi_{22}
	\end{bmatrix}\begin{bmatrix}
	I & X_{+}-BU \\
	0 & -X_{-} \\
	\end{bmatrix}^{\top}\label{equ:N}.
\end{align}
Now, we examine that $M$ and $N$ satisfy the kernel assumption in Lemma \ref{lem:matrixS}. Note that $M_{22}\leq 0$ and $N_{22} \leq 0$ trivially hold since $Z= Y-X^{\top}X > 0$ and $\Phi_{22} <0$. Clearly, $\text{ker}(N_{12}) = \text{ker} (\Phi_{12} + (X_{+}-BU)\Phi_{22})X_{-}^{\top}$, and $\text{ker}(N_{22}) = \text{ker}(X_{-}^{\top})$ as $\Phi_{22}<0$. Thus, $\text{ker}(N_{22}) \subset \text{ker}(N_{12})$. The generalized Slater condition implies that the set $\Sigma$ has at least an interior point, which is a mild assumption in our problem. Then, we have the following main result under the Slater condition.
\begin{theorem}\label{theorem:Lmi}
	Assume that there exists some matrix $A$ such that
	$
	\begin{bmatrix}
	I \\
	A^{\top} 
	\end{bmatrix}^{\top}
	N
	\begin{bmatrix}
	I \\
	A^{\top} 
	\end{bmatrix} > 0.
	$
	Then, the system is stabilizable via logarithmically quantized linear feedback with density $\rho = (1-\delta)/(1+\delta)$ for all $A \in \Sigma$ if and only if there exists $Y>0$, $X$ and scalars $\alpha \geq 0, \beta >0, \delta > 0$ such that the following LMI holds
	\begin{equation}\label{equ:LMI}
	\hspace{-0.085cm}\begin{bmatrix}
	Y-\delta^2BB^{\top}-\beta I & 0 & BX & 0 \\
	0 & 0& Y & 0\\
	X^{\top}B^{\top} & Y^{\top} & Y & X^{\top} \\
	0 & 0 & X & I
	\end{bmatrix}
	-\alpha 
	\begin{bmatrix}
	I & X_{U} \\
	0 & -X_{-} \\
	0 & 0\\
	0 & 0
	\end{bmatrix}\begin{bmatrix}
	\Phi_{11} & \Phi_{12} \\
	\Phi_{12}^{\top} & \Phi_{22}
	\end{bmatrix}\begin{bmatrix}
	I & X_{U} \\
	0 & -X_{-} \\
	0 & 0\\
	0 & 0
	\end{bmatrix}^{\top} \geq 0,
	~~\begin{bmatrix}
	Y & X^{\top} \\
	X & I
	\end{bmatrix} > 0.
	\end{equation}
	Moreover, if (\ref{equ:LMI}) is feasible for some $Y$ and $X$, then a stabilizing controller is given by $u(k) = f(v(k)) $ with quantization density $\rho= (1-\delta)/(1+\delta)$ and $v(k)= Kx(k)$ with $K = XY^{-1}$. 
\end{theorem}
\begin{proof}
	To prove the ``if'' statement, suppose that (\ref{equ:LMI}) is feasible. Now, we calculate the Schur complement of the first LMI in (\ref{equ:LMI}) concerning $I$ and obtain
	$$
	\begin{bmatrix}
	Y-\delta^2BB^{\top}-\beta I & 0 & BX  \\
	0 & 0& Y \\
	X^{\top}B^{\top} & Y^{\top} & Z \\
	\end{bmatrix}
	-\alpha 
	\begin{bmatrix}
	I & X_{U} \\
	0 & -X_{-} \\
	0 & 0\\
	\end{bmatrix}\begin{bmatrix}
	\Phi_{11} & \Phi_{12} \\
	\Phi_{12}^{\top} & \Phi_{22}
	\end{bmatrix}\begin{bmatrix}
	I & X_{U} \\
	0 & -X_{-} \\
	0 & 0\\
	\end{bmatrix}^{\top} \geq 0,
	~~\begin{bmatrix}
	Y & X^{\top} \\
	X & I
	\end{bmatrix} > 0.
	$$
	Then, we again compute the Schur complement with respect to $Z$ to yield $M-\alpha N \geq \begin{bmatrix}
	\beta I & 0 \\
	0 & 0
	\end{bmatrix}$, where we have used the fact that $Z = Y - X^{\top}X >0$ by the second inequality of (\ref{equ:LMI}). Thus, by the matrix S-lemma in Lemma \ref{lem:matrixS}, the inequality (\ref{equ:performance}) holds for all $A\in \Sigma$. According to our change of variables, a stabilizing controller $K$ can be solved by $K = XY^{-1}$.
	
	To prove the ``only if'' statement, suppose that the system is stabilizable with quantization density $\rho$ for all $A \in \Sigma$. 
	By assumption, the inequality (\ref{equ:performance}) holds for all $A\in \Sigma$ with $Z = Y- X^{\top}X > 0$. Thus, the first inequality in (\ref{equ:LMI}) follows from Lemma \ref{lem:matrixS} and some Schur complement. The second inequality in  (\ref{equ:LMI}) holds since $Z = Y- X^{\top}X > 0$. The proof is now completed.
\end{proof}

It can be clearly observed from the first LMI in (\ref{equ:LMI}) that as $\delta$ increases, the system becomes harder to stabilize. Hence, there exists a coarsest quantization density such that the system cannot be stabilized for any lower density.

\section{Coarsest quantization density}
In this section, we derive a tight lower bound of the quantization density $\rho$ for the system to be stabilizable, provided with the data matrices $\{X_{-}, X_{+}, U\}$.

Consider the following min-max $\mathcal{H}_{\infty}$ optimization problem
\begin{equation}\label{prob:minmax}
\min \limits_{K} \max \limits_{A \in \Sigma} ~~\|G_{A,K}(z) \|_{\infty}.
\end{equation}
Let $\mathcal{K}$ be the set of optimal controllers, $\mathcal{S}$ be the set of optimal state matrices, and $\gamma^*$ be its optimal value to (\ref{prob:minmax}). Then, we show that $\delta =(1-\rho)/(1+\rho)$ is in fact upper bounded by $1/\gamma^*$. 

\begin{theorem}
	Suppose that (\ref{prob:minmax}) is feasible. Then, there exists a stabilizing logarithmically quantized controller with density $\rho=(1-\delta)/(1+\delta)$ for all systems with $A \in \Sigma$ if and only if $\delta < 1/ \gamma^*$.
\end{theorem}

\begin{proof}
	To prove the ``if'' statement, suppose that $\delta < 1/\gamma^*$. Let $K^* \in \mathcal{K}$ be an optimal controller to (\ref{prob:minmax}). Then, for any $A \in \Sigma$, it follows that 
	$$\|G_{A,K^*}(z)\|_{\infty} \leq  \max \limits_{A\in \Sigma} \|G_{A,K^*}(z)\|_{\infty} = \gamma^*.$$
	Multiplying by $\delta$ in both sides of the above inequality yields that $\delta\|G_{A,K^*}(z)\|_{\infty}< 1$ by the assumption $\delta < 1/\gamma^*$. Hence, by (\ref{equ:h_inf}), the feedback gain $K^*$ is able to stabilize (\ref{equ:sys}) with quantization density $\rho = (1-\delta)/(1+\delta)$.
	
	To prove the ``only if'' statement, suppose that a controller $K$ with quantization density $\rho$ stabilizes all systems with $A\in \Sigma$. Let $A^* \in \mathcal{S}$ be an optimal system matrix to  (\ref{prob:minmax}). Then, it follows that
	$$
	 \delta \gamma^* = \delta \min \limits_{K} \|G_{A^*,K}(z)\|_{\infty}  < \delta \|G_{A^*,K}(z)\|_{\infty} <1.
	$$
	Therefore, we must have that $\delta < 1/\gamma^*$. The proof is now completed.
\end{proof}
Note that (\ref{prob:minmax}) may not have a solution if $\Sigma$ is unbounded. In fact, the boundedness of $\Sigma$ is a necessary condition for (\ref{prob:minmax}) to be feasible, as to be shown in the next section. Now, we propose an efficient SDP to solve (\ref{prob:minmax}). Clearly, (\ref{prob:minmax}) is equivalent to 
$$
\min \limits_{K,\delta>0} 1/\delta, ~~~~~ \text{subject to} ~~\|G_{A,K}(z) \|_{\infty} < 1/\delta, ~~ \forall A \in \Sigma.
$$
The $\mathcal{H}_{\infty}$ norm constraint in the above minimization problem can be formulated as an LMI by virtue of Theorem \ref{theorem:Lmi}, which leads to the following maximization problem 
\begin{equation}\label{equ:SDP}
\max \limits_{Y,X,\alpha,\beta,\delta} \delta, ~~~~~\text{subject to} ~~(\ref{equ:LMI}) ~~\text{and}~~\delta>0, Y>0, \alpha \geq 0, \beta >0,
\end{equation}
which is an SDP and can be efficiently solved by modern solvers e.g., CVX~\citep{cvx}.

\section{Necessary condition for quantized stabilization}

So far, we have established both sufficient and necessary conditions for quantized stabilization and have proposed an SDP to solve the coarsest quantization density, all based on the LMI in (\ref{equ:LMI}). To provide more profound insights, we derive a necessary and explicit rank condition on the data matrices by leveraging classical quantized control results.

Consider the following max-min problem that has interchanged the order of minimization and maximization of (\ref{prob:minmax})
\begin{equation}\label{prob:maxmin}
\max \limits_{A\in \Sigma} \min \limits_{K} \|G_{A,K}(z) \|_{\infty}.
\end{equation}
According to \cite{fu2005sector}, the inner minimization can be expressed as
$
\min \limits_{K} \|G_{A,K}(z) \|_{\infty} =  \prod_{i}|\lambda_i| 
$,
where $\lambda_i$ denotes the $i$-th unstable eigenvalue of $A$. Thus, (\ref{prob:maxmin}) is equivalent to a maximization problem
$$
\max \limits_{A \in \Sigma} \prod_{i}|\lambda_i|.
$$

By the well-known max-min inequality~\citep{boyd2004convex}, it follows that
$$
\max \limits_{A \in \Sigma} \prod_{i}|\lambda_i|  \leq  \min \limits_{K} \max \limits_{A \in \Sigma} \|G_{A,K}(z) \|_{\infty} = \gamma^*.
$$
Hence, the maximum $\max \limits_{A \in \Sigma} \prod_{i}|\lambda_i|$ must be bounded such that $\gamma^*$ is finite. Since $\Sigma$ is only related to the data, we then ask the following question: under which conditions on the data matrices, the eigenvalues of $A$ are bounded for all $A \in \Sigma$? 

A straightforward sufficient condition is that $\Sigma$ is bounded, which is equivalent to $\text{rank}(X_{-}) = n$. However, it is non-trivial to prove that the rank condition is also necessary for the eigenvalue to be bounded. The following example shows that when $\Sigma$ is unbounded, there might still exist a nilpotent matrix $A \in \Sigma$ with unbounded elements yet zero eigenvalues.

\begin{example}
	Consider a nilpotent matrix $
	A = \begin{bmatrix}
		0 & k \\
		0 & 0
	\end{bmatrix}
	$ and input matrix $B = I$. The data matrices are given by $X_{-} = \begin{bmatrix}
	1 & 1\\
	0 & 0
	\end{bmatrix}$, $U = I$, $X_{+} = I$. It can be easily verified that $X_{+} = AX_{-}+ BU$. Set $\Phi_{11} = I, \Phi_{12}=0, \Phi_{22}=-I$. Then, it is easy to verify that $A \in \Sigma$ for $k \in \mathbb{R}$. Hence, $A$ can be unbounded but its eigenvalues are bounded.
\end{example}

Nevertheless, we show by contradiction that the boundedness of $\Sigma$ is both sufficient and necessary for the eigenvalues of $A$ to be bounded for all $A \in \Sigma$, which leads to the following result.

\begin{theorem}
	If there exists a stabilizing logarithmically quantized controller, then 
	$
	\text{rank}(X_{-}) = n.
	$
\end{theorem}
\begin{proof}
	We prove by contradiction. Suppose that $\text{rank}(X_{-}) = r <  n$. Next, we construct a matrix $\bar{A} \in \Sigma$ which has unbounded eigenvalues. First, we define an invertible linear transformation $E = \begin{bmatrix}
	E_{11} & E_{12} \\
	E_{21} & E_{22}
	\end{bmatrix} \in \mathbb{R}^{n \times n} $ such that $EX_{-} = \begin{bmatrix}
	X_{r} \\
	\hline 0 
	\end{bmatrix}
	$ with $X_r \in \mathbb{R}^{r \times T}$, and a diagonal matrix $\Lambda = \text{diag}(0,\dots,0,1,\dots,1)  \in \mathbb{R}^{n \times n} $ with its first $r$ diagonal elements being zero. 
	
	Let $A_0 \in \Sigma$, and $\bar{A} = A_0 + k \Lambda E$. Then, we prove that $\bar{A} \in \Sigma$. To see this,  note that $\Lambda E X_{-} = 0$ by definition, and $\Sigma$ can be written as 
	$$
	\begin{bmatrix}
	I & X_{U}+AX_{-} 
	\end{bmatrix}\begin{bmatrix}
	\Phi_{11} & \Phi_{12} \\
	\Phi_{12}^{\top} & \Phi_{22}
	\end{bmatrix}\begin{bmatrix}
	I & X_{U}+AX_{-} 
	\end{bmatrix}^{\top} \geq 0.
	$$
	
	Next, we show that $\Lambda E$ has a non-zero eigenvalue. By the definition, $E_{22}$ has full rank. Let $\lambda_{E}$ be an eigenvalue of $E_{22}$ and $x_{E}$ be the associated eigenvector. Denote $x = \begin{bmatrix}
	0 \\
	\hline x_E
	\end{bmatrix} \in \mathbb{R}^{n}$. Then, it is straightforward to verify that $\Lambda E x = \lambda_E x$, which implies that $\lambda_E \neq 0$ is an eigenvalue of $\Lambda E$.
	
	Hence, by limiting $k$ to infinity, we conclude that $\bar{A}$ has unbounded eigenvalue. Thus, we have that $\max \limits_{A \in \Sigma} \prod_{i}|\lambda_i| = \infty$, which leads to a contradiction.
\end{proof}
Though the condition $\text{rank}(X_{-}) = n$ is only necessary, it provides a simple approach to examine the existence of a solution to (\ref{equ:LMI}). Moreover, it also holds for zero quantization density, which corresponds to general linear systems without quantization. Hence, our results extend \cite{van2020noisy} in that we provide a necessary and explicit rank condition for quadratic stabilization.

\section{Numerical examples}
\begin{figure}[t]
	\centering
	\includegraphics[height=50mm]{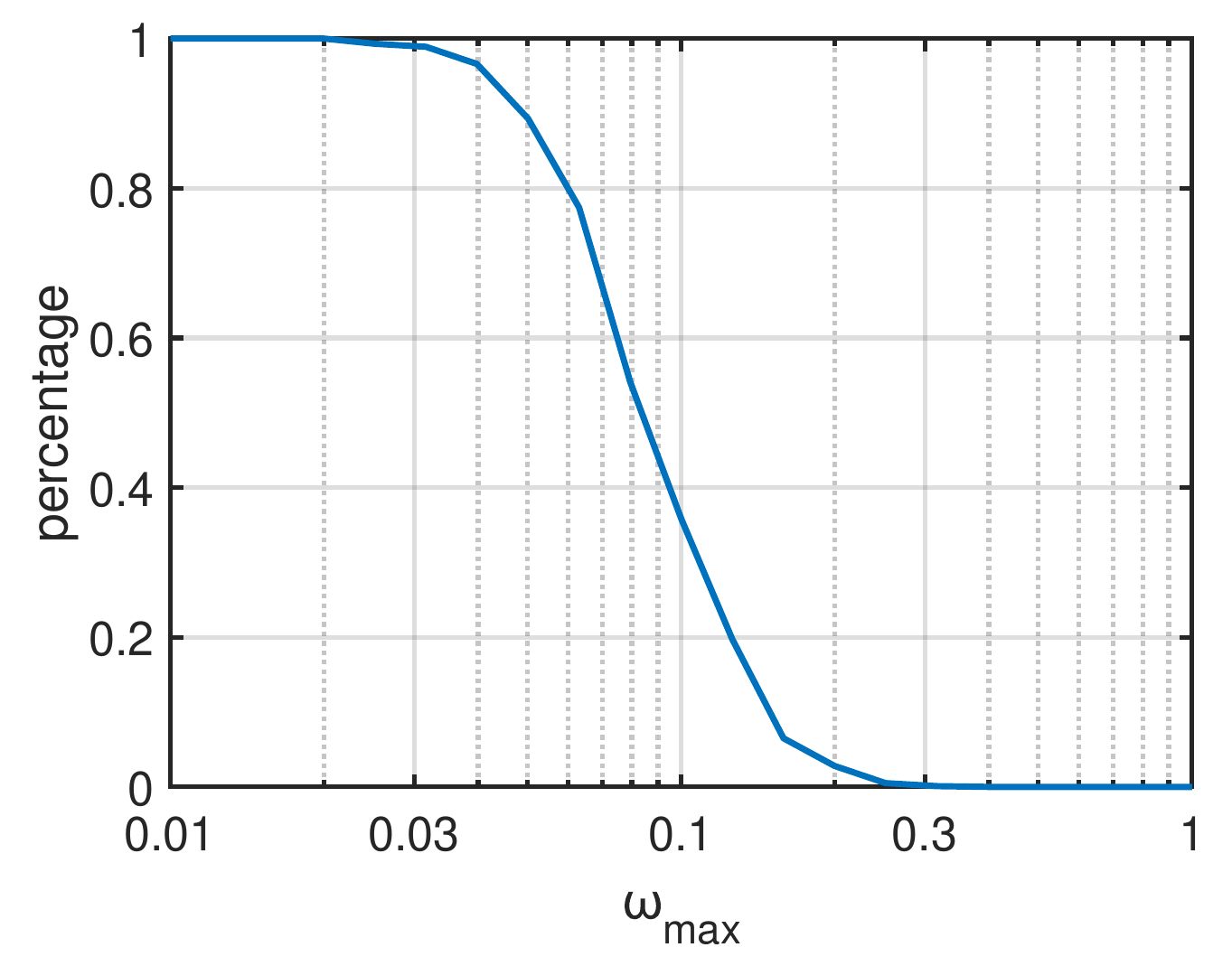}
	\includegraphics[height=50mm]{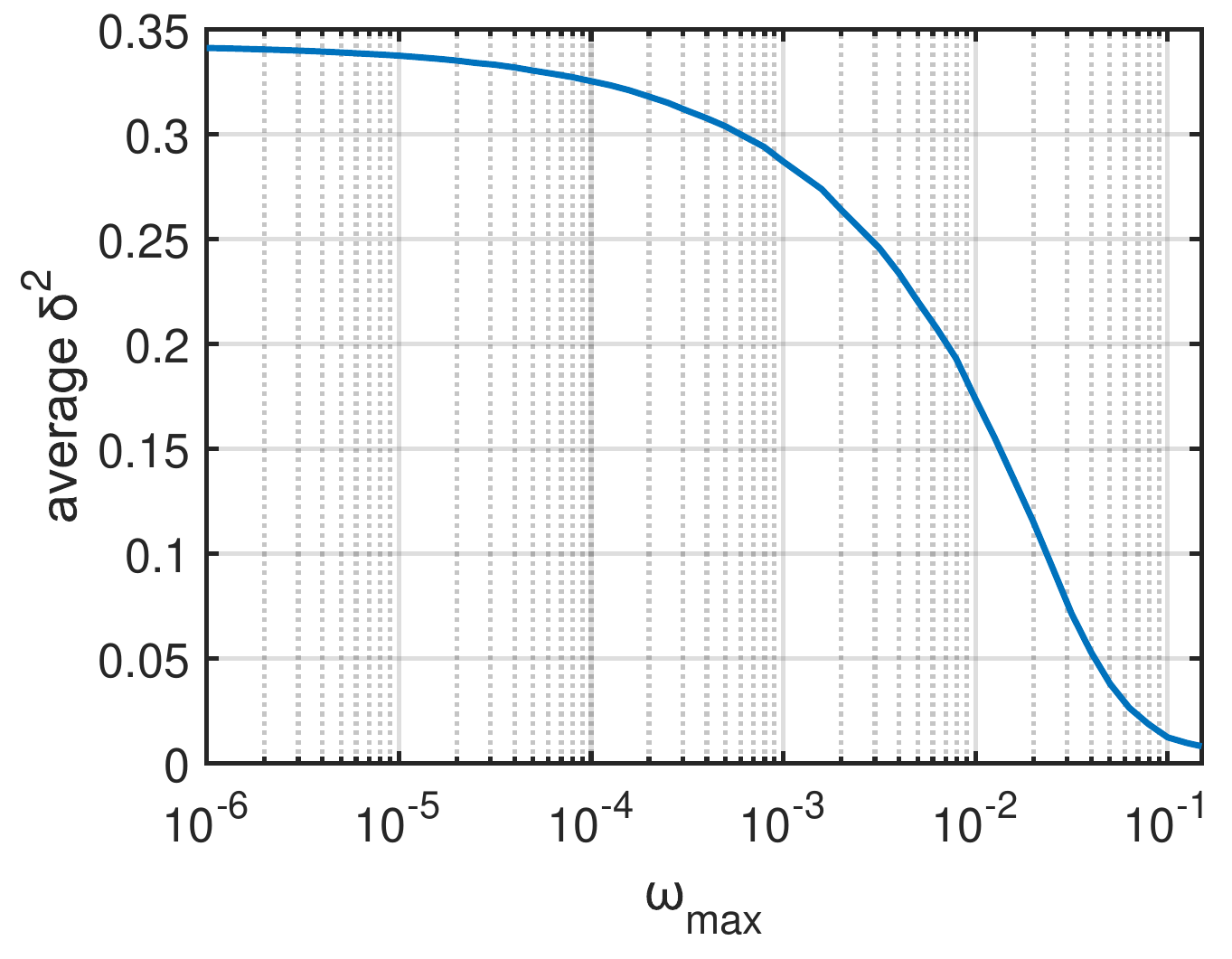}
	\caption{Left: the percentage of data sets for which (\ref{equ:LMI}) has a solution; Right: the mean of $\delta^2$ which is obtained by solving (\ref{equ:SDP}). The x-axis is set to be logarithmic for better exposition.}
	\label{pic:noise}
\end{figure}
In this section, we perform numerical examples to verify our theoretical results. The simulation is carried out using MATLAB 2020b on a laptop with a 2.8GHz CPU. The code is provided in \url{https://github.com/lixc21/Data-driven-Quantized-Control}.

We randomly generate an open-loop unstable linear dynamical model
$$
A=\begin{bmatrix}
-0.192 & -0.936 & -0.814 \\
-0.918 & 0.729 & -0.724 \\
-0.412 & 0.735 & -0.516
\end{bmatrix},~
B=\begin{bmatrix}
-0.554 \\ 0.735 \\ 0.528
\end{bmatrix}
$$ 
by sampling its elements uniformly from the interval $[-1,1]$. The eigenvalues of $A$ are $1.2910$, $-1.3228$, $0.0528$. The time horizon of the trajectory is set to $T = 20$. We assume that the noise is independently sampled from a uniform distribution on a three-dimensional ball $\{w\in \mathbb{R}^3|~\|w\|_2^2\leq \omega_{\text{max}}\}$ with noise level $\omega_{\text{max}}$. The prior bound in (\ref{assumption}) on the noise is then given by $\Phi_{11}=T\omega_{\text{max}}I$, $\Phi_{12}=0$, $\Phi_{22}=-I$. It can be easily verified that the sampled noises satisfy the quadratic bound. The initial state and the control input at each time instant are sampled from the standard normal distribution. 

First, we investigate the effects of noise level on quantized stabilization. We let $\omega_{\text{max}}$ vary from $0$ to $1$. For each noise level, we generate $1000$ independent data sets to solve the LMI in (\ref{equ:LMI}) and plot the percentage for (\ref{equ:LMI}) to have a feasible solution, i.e., a stabilizing controller exists. We check that the generalized Slater condition required by Theorem \ref{theorem:Lmi} holds for all the data sets by verifying that $N$ in (\ref{equ:N}) has three positive eigenvalues. The result is displayed on the left of Fig. \ref{pic:noise}. Consistent with intuition, the percentage decreases as the noise level grows. This is because the uncertainty set of $A$ becomes larger and it is harder to find a quantized stabilizing controller. We obtain similar results for the coarsest quantization density in the right of Fig. \ref{pic:noise}.


\begin{figure}[t]
	\centering
	\includegraphics[height=50mm]{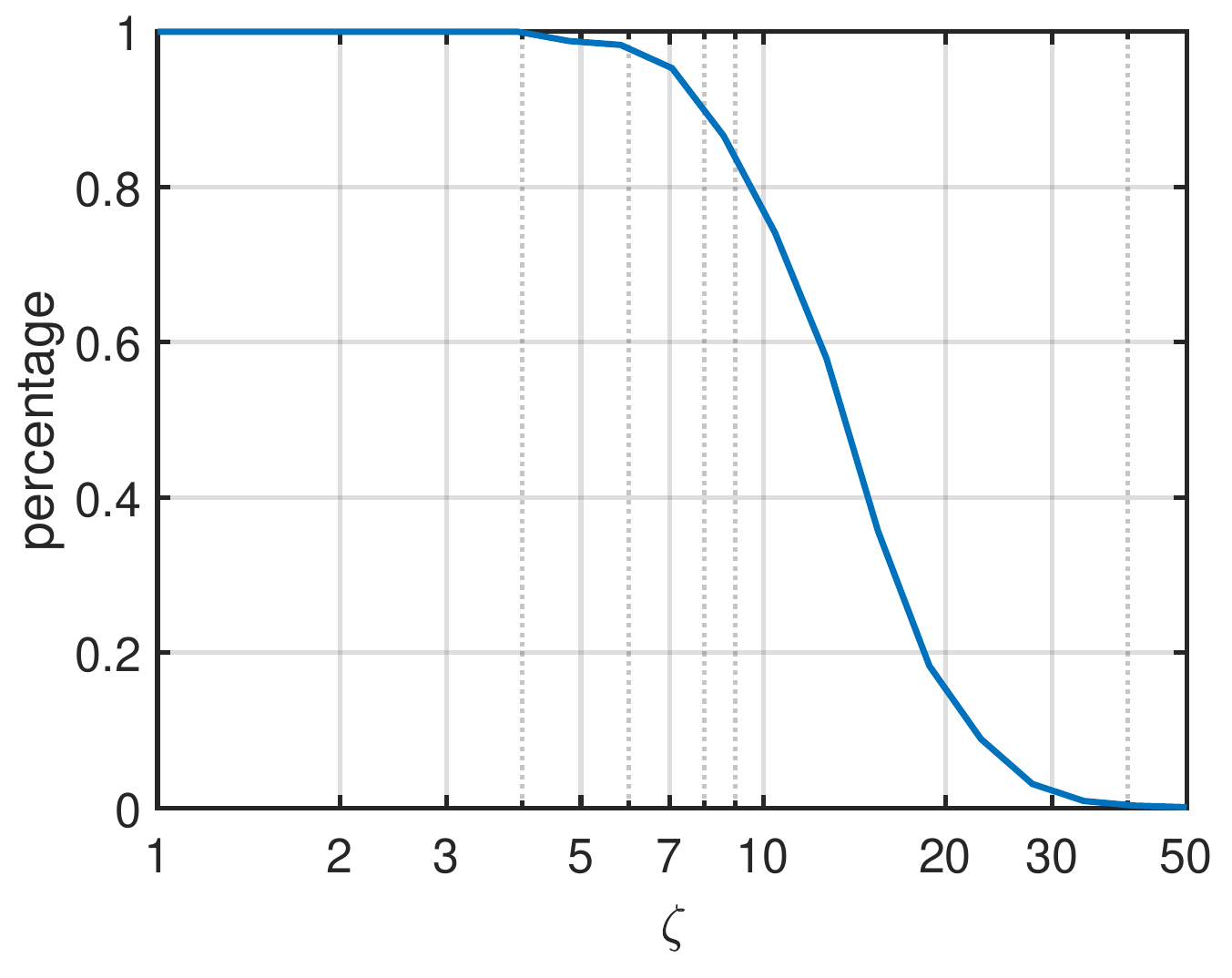}
	\includegraphics[height=50mm]{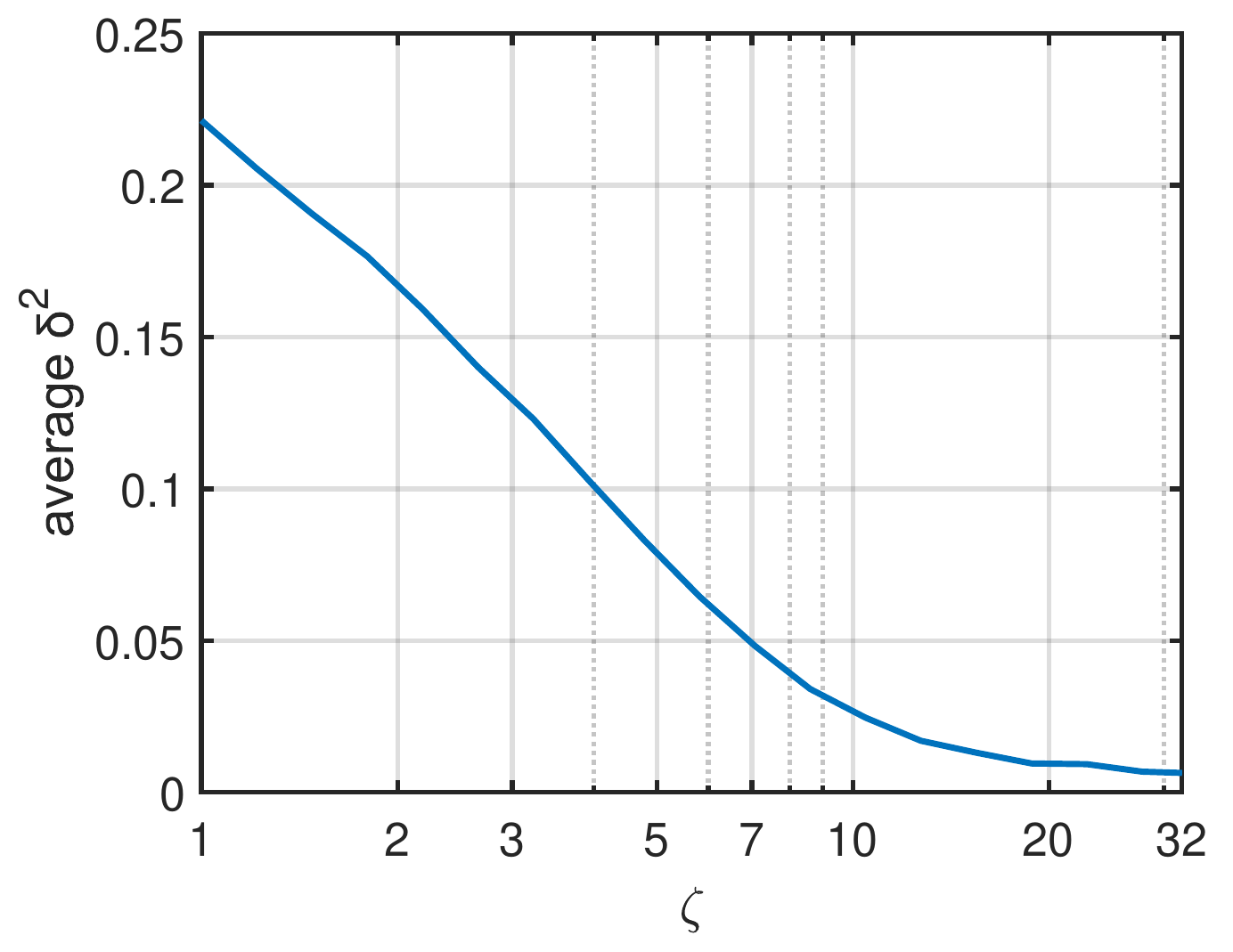}
	\caption{Left: the percentage of data sets for which (\ref{equ:LMI}) has a solution; Right: the mean of $\delta^2$ which is obtained by solving (\ref{equ:SDP}). The x-axis is set to be logarithmic.}
	\label{pic:bound}
\end{figure}

Then, we note that the prior bound on the noise in (\ref{assumption}) can also have a significant impact on the performance. To see this, we let the prior bound be larger than the true bound by setting $\Phi_{11}=\zeta \times T\omega_{\text{max}}I$, $\Phi_{12}=0$, $\Phi_{22}=-I$, where $\omega_{\text{max}}=0.005$ and the parameter $\zeta$ varies from 1 to 50. For each $\zeta$, we generate 1000 data sets independently and compute the percentage for (\ref{equ:LMI}) to have a solution. For the cases that a stabilizing controller exists, we further solve the coarsest quantization density by (\ref{equ:SDP}). The generalized Slater condition is also verified to hold for all data sets. The results are shown in Fig. \ref{pic:bound}. Clearly, if we have a more accurate prior bound, we are more likely to stabilize the system via a lower quantization density.

\section{Conclusion}

In this paper, we have addressed the stabilization problem using logarithmically quantized feedback for partially unknown linear systems with noisy inputs. Particularly, we have provided a sufficient and necessary condition for quantized stabilization in the form of an LMI. Moreover, we have solved the coarsest quantization density as well as a stabilizing feedback gain via an SDP. 

We believe that our results raise a range of interesting directions. For example, one may extend our results to the stabilization of general MIMO systems. Instead of solving a stabilizing controller, we can also study $\mathcal{H}_2$ and $\mathcal{H}_{\infty}$ performance guarantees of the closed-loop systems. For online control with quantized feedback, one can design an adaptive quantizer with a varying density as more data is collected per round. 

\section{Acknowledgement}
We gratefully acknowledge support from the National Natural Science Foundation of China under Grant no. 62033006. We also would like to thank the anonymous reviewers for their helpful suggestions.

\bibliography{mybibfile}

\end{document}